\documentclass{article}
\usepackage{comment}
\usepackage{amsmath,amssymb,amsthm, amsfonts}
\usepackage{url}
\usepackage{tikz}
\usetikzlibrary{decorations}
\usetikzlibrary{decorations.pathmorphing}
\usetikzlibrary{decorations.pathreplacing}
\usetikzlibrary{decorations.markings}
\usetikzlibrary{decorations.text}
\usetikzlibrary{arrows.meta}
\usetikzlibrary{calc}
\usetikzlibrary[arrows,positioning, shapes, fit]

\newcommand{\Z}{\operatorname{Z}}
\newcommand{\Zo}{\operatorname{Z^-}}
\newcommand{\F}{\operatorname{F}}
\newcommand{\Fo}{\operatorname{F^-}}

\usepackage{xcolor}
\date{}
\usepackage{authblk}
\newtheorem{thm}{Theorem}
\newtheorem{lem}[thm]{Lemma}
\newtheorem{dfn}[thm]{Definition}
\newtheorem{cor}[thm]{Corollary}

\newtheorem{obs}[thm]{Observation}
\newtheorem{prop}[thm]{Proposition}
\title{The zero forcing span of a graph}
\author{Bonnie Jacob}
\affil{National Technical Institute for the Deaf,  Rochester Institute of Technology \\
Rochester, NY, USA, bcjntm@rit.edu}
\begin{document}

\maketitle              

\begin{abstract}
In zero forcing, the focus is typically on finding the minimum cardinality of any zero forcing set in the graph; however, the number of cardinalities between $0$ and the number of vertices in the graph for which there are both zero forcing sets and sets that fail to be zero forcing sets is not well known.  In this paper, we introduce the zero forcing span of a graph, which is the number of distinct cardinalities for which there are sets that are zero forcing sets and sets that are not.  We introduce the span within the context of standard zero forcing and skew zero forcing as well as for standard zero forcing on directed graphs.  We characterize graphs with high span and low span of each type, and also investigate graphs with special zero forcing polynomials.  

\end{abstract}


\section{Introduction}

Throughout this paper, we use $G$ to denote a finite, simple graph on $n=|V(G)|$ vertices with edge set $E(G)$.  We use $D$ to denote a directed graph, or digraph, with vertex set $V(D)$ and arc set $E(D)$.  Like in our undirected graphs, hereafter simply ``graphs,'' we do not allow loops or multiple arcs in our digraphs, though between a pair of vertices there may be an arc in each direction (from vertex $u$ to vertex $v$ and from vertex $v$ to vertex $u$, for example).  When the graph is understood, we use $V$ in place of $V(D)$ or $V(G)$.  Unless otherwise stated, we use $n=|V(G)|$ and call $|V(G)|$ the \emph{order} of the graph.  

For any $v \in V(G)$, the \emph{open neighborhood} of $v$ denoted $N(v)$ is the set of vertices adjacent to $v$.  In a digraph, the \emph{open in-neighborhood} of $v$ denoted $N^-(v)$ is the set of vertices from which there is an arc to $v$, and the \emph{open out-neighborhood} denoted $N^+(v)$ the set of vertices to which there is an arc from $v$.  A \emph{neighbor} of $v$ is a vertex in the open neighborhood of $v$, with analogous definitions for in- and out-neighbor.  

Zero forcing is a process that consists of designating a subset $S \subseteq V(G)$ as blue, and the remaining vertices as white.  A color change rule of varying forms is then applied.  If repeated applications of the color change rule results in all vertices eventually turning blue, the original set is called a \emph{zero forcing set}.  

In this paper, we use three different color change rules, and define those here.  

\begin{enumerate}
\item The \emph{standard color change rule} (standard zero forcing): if any blue vertex has exactly one white neighbor, then the white neighbor becomes blue.
\item The \emph{skew color change rule} (skew zero forcing): if any  vertex (blue or white) has exactly one white neighbor, then the white neighbor becomes blue.  
\item The \emph{(standard) digraph color change rule} (standard zero forcing on digraphs): if any blue vertex  has exactly one white out-neighbor, then the white out-neighbor becomes blue.
\end{enumerate}

Under each of these rules, the minimum cardinality of any starting set of blue vertices that eventually results in the entire graph turning blue is called  the \emph{zero forcing number}, first formally introduced in \cite{aim2008zero} (the \emph{skew zero forcing number} introduced in \cite{ima2010minimum}, or the \emph{digraph zero forcing number} introduced in \cite{berliner2013minimum} respectively) and is denoted $\Z(G)$ (or $\Zo(G)$, or $\Z(D)$).  The maximum cardinality of any starting set of blue vertices that never results in the entire graph turning blue  is called the \emph{failed zero forcing number}  introduced in \cite{fetcie2015failed} (the \emph{failed skew zero forcing number} introduced in \cite{ansill2016failed}, or the \emph{digraph failed zero forcing number} introduced in \cite{adams2021failed} respectively) and is denoted $\F(G)$ (or $\Fo(G)$, or $\F(D)$).   

Given a graph $G$, the \emph{zero forcing span} $\lambda(G)$ is the number of distinct cardinalities, $k_1, k_2, \ldots k_{\lambda(G)}$ such that for each $i$, $1 \leq i \leq \lambda(G)$,  there exist sets $Z, F \subseteq V(G)$ with $|Z|=|F|=k_i$, where $Z$ is a zero forcing set and $F$ is not. We define $\lambda^-(G)$ analogously, but where $Z$ is a skew zero forcing set and $F$ is not, and $\lambda(D)$ as well, but where $D$ is a digraph and the digraph color change rule is applied.

Note that $\lambda(G)$ denotes how much bigger than $\Z(G)$  a set must be  to guarantee that it must be a zero forcing set, without regard to which vertices are in the set.

\section{Motivation}

The concept of zero forcing span has connections with several other problems in the literature.  First, we look at the connection to linear algebra.  With $G$ we associate a set of symmetric matrices denoted $\mathcal{S}(G)$.  Number the vertices $1, 2, \ldots, n$ and define  $\mathcal{S}(G)$ as follows.  
$$\mathcal{S}(G) = \left\{ A \in \mathbb{R}^{n \times n} \ : \ A^T = A \ \mbox{ and for } i \neq j, a_{ij} \neq 0 \mbox{ if and only if } ij \notin E(G) \right\}.$$
Note that the diagonal is unconstrained.  

By applying \cite[Proposition 2.3]{aim2008zero}, we find that the zero forcing span $\lambda(G)$ is  the number of values $k$ such that both of the following conditions hold.  
\begin{enumerate}
\item There exists a set $S \subseteq \{1, 2, \ldots, |V(G)|\}$ with $|S|=k$ such that for any $A \in \mathcal{S}(G)$, if $v \in \ker(A)$ with $v_i=0$ for all $i \in S$, then $v=0$, but
\item there exists $A \in \mathcal{S}(G)$ with $v \in \ker(A)$, $v \neq 0$ and some $k$ entries of $v$ that are all $0$.  
\end{enumerate}

We can also relate the zero forcing span to the \emph{zero forcing polynomial} $\mathcal{Z}(G;x)=\sum_{i=1}^n z(G;i) x^i$, introduced in \cite{boyer2019zero}, where $z(G;i)$ is the number of zero forcing sets of cardinality $i$.  Then $\lambda(G)$ gives the number of terms where $0< z(G;i) <{ n \choose i }$, that is, the number of terms in the polynomial that are not simply $0$  or ${n \choose i}$ (the minimum or maximum possible for each coefficient in the polynomial).  

While not yet explicitly defined in the literature, we can define analogous zero forcing polynomials for variants of zero forcing.
Let the \emph{skew zero forcing polynomial} be $\mathcal{Z^-}(G;x)=\sum_{i=0}^n z^-(G;i) x^i$ where $z^-(G;i)$ is the number of skew zero forcing sets of cardinality $i$, and the \emph{directed zero forcing polynomial} be $\mathcal{Z^D}(D;x)=\sum_{i=1}^n z^D(D;i) x^i$ where $z^D(D;i)$ is the number of zero forcing sets of cardinality $i$ in a digraph $D$. Note that, unlike standard zero forcing, the skew zero forcing polynomial may have $z(G;0) > 0$.

\section{Results}

\subsection{Extreme values of $\lambda(G)$} 

We  note the following formula for $\lambda(G)$ in terms of $\F(G)$ and $\Z(G)$.  The equivalent statement holds for each type of zero forcing.  

\begin{obs}
$\lambda(G)=\F(G) - \Z(G) +1$.  \label{obs:formula}
\end{obs}

Recall that  $\Z(G)\geq 1$ for any undirected graph $G$.  Also recall that $\F(G) \leq n-1$ and that $\F(G) \geq \Z(G)-1$.  The same statements hold for a digraph $D$, which give us the following trivial bounds.  

\begin{obs}
For a graph $G$ and digraph $D$, $0 \leq \lambda(G) \leq n-1$ and $0 \leq \lambda(D) \leq n-1$.
\end{obs}

In the skew case, there exist graphs for which $\Zo(G)=0$, implying the following observation.

\begin{obs}
$0 \leq \lambda^-(G) \leq n$ \label{obs:skew}
\end{obs}

\subsubsection{Characterizations of graphs with zero forcing span of $0$}

\begin{lem}
The following are equivalent:
\begin{enumerate}
\item $\lambda(G)=0$
\item $\F(G) < \Z(G)$
\item $\Z(G) = \F(G) +1$
\end{enumerate}
The same equivalence holds if the graph $G$ is replaced by a digraph $D$.  \label{lem:equivstandard}
\end{lem}
\begin{proof}
By Observation \ref{obs:formula}, $\lambda(G)=0$ if and only if $\Z(G) = \F(G)+1$.  Since $\F(G) \geq \Z(G)-1$ (as noted in \cite{fetcie2015failed} or by observing that any set of vertices must be a failed zero forcing set or a zero forcing set), the equivalence holds.  \qed
\end{proof}

For skew zero forcing, since $\Fo(G)$ is not always defined, we have the following list of equivalences.  The proof is identical to that of Lemma \ref{lem:equivstandard} but with the addition of the possibility that we may have $\Zo(G)=0$, which is equivalent to $\Fo(G)$ being undefined.  
\begin{lem}
The following are equivalent:
\begin{enumerate}
\item $\lambda^-(G)=0$
\item $\Fo(G) < \Zo(G)$ or $\Zo(G)=0$.  
\item $\Zo(G) = \Fo(G) +1$ or $\Fo(G)$ is undefined.  
\end{enumerate}  \label{lem:equivskew}
\end{lem}

In \cite{fetcie2015failed}, it was established that $\F(G)<\Z(G)$ if and only if $G=K_n$ or $G=\overline{K_n}$, leading to the following characterization of graphs with $\lambda(G)=0$.    

\begin{thm}
$\lambda(G)=0$ if and only if  $G=K_n$ or $G=\overline{K_n}$ \label{prop:zero}
\end{thm}

\begin{proof}
We have that $\F(G)=\Z(G)-1$ if and only  if $G=K_n$ or $G=\overline{K_n}$ \cite{fetcie2015failed}.  By Lemma \ref{lem:equivstandard}, the result follows. \qed 
\end{proof}

For skew zero forcing, we  need a few definitions to characterize graphs with $\lambda^-(G)=0$.  In \cite{ansill2016failed}, we established that $\Fo(G)<\Zo(G)$ if and only if $G$ is an odd cycle or nonempty set of cycles intersecting in a single vertex, or a doubly extended bouquet-dipole, pictured in Figure \ref{fig:skewzero} on the left.   

\begin{dfn}
We call a graph $G$ a \emph{doubly extended bouquet-dipole} if it consists of vertices $u$ and $v$ that are each on a nonempty set of odd cycles, where all other vertices on the cycles have degree two, and $u, v$ are joined by a path of even order that alternates between single even order paths whose internal vertices all have degree two, and multiple even order paths whose internal vertices all have degree two.  
\end{dfn} 
\tikzset{
  c/.style={every coordinate/.try}
}

\begin{figure}[h!]
\begin{center}
\begin{tikzpicture}[auto, scale=0.825]
\tikzstyle{vertex}=[draw, circle, inner sep=0.8mm]
\node (v1) at (360/5: 0.5) [vertex] {};
\node (v2) at (2*360/5:0.5) [vertex] {};
\node (v3) at (3*360/5: 0.5) [vertex] {};
\node (v4) at (4*360/5: 0.5) [vertex] {};
\node (v5) at (5*360/5:0.5) [vertex, label=left:$u$] {};
\node (v6) at (0.5, -0.6) [vertex] {};
\node (v7) at (0.9, -0.4) [vertex] {};
\node (v8) at (0.5, 0.6) [vertex] {};
\node (v9) at (0.9, 0.4) [vertex] {};

\node (v10) at (1.6,0) [vertex]{};
\node (v11) at (2.4,0) [vertex]{};
\node (v12) at (3.2,0) [vertex]{};

\node (v13) at (4.0,0) [vertex]{};
\node (v14) at (4.8,0) [vertex]{};
\node (v15) at (5.6,0) [vertex]{};

\node (v16) at (3.5,0.6) [vertex]{};
\node (v17) at (4.1, 0.6) [vertex]{};
\node (v18) at (4.7,0.6) [vertex]{};
\node (v19) at (5.3,0.6) [vertex]{};

\node (v20) at (4.0,-0.6) [vertex]{};
\node (v21) at (4.8,-0.6) [vertex]{};

\node (v22) at (6.4, 0) [vertex]{};
\node (v23) at (7.2, 0) [vertex]{};
\node (v24) at (8.0, 0) [vertex, label=right:$v$]{};

 \node (v25) at ([shift={(8.6,0)}]231: 0.6) [vertex] {};
\node (v26) at ([shift={(8.6,0)}]283: 0.6) [vertex] {};
 \node (v27) at ([shift={(8.6,0)}]334: 0.6) [vertex] {};
\node (v28) at ([shift={(8.6,0)}]26: 0.6) [vertex] {};
 \node (v29) at ([shift={(8.6,0)}]77: 0.6) [vertex] {};
\node (v30) at ([shift={(8.6,0)}]129: 0.6) [vertex] {};

\foreach[evaluate=\y using int(\x-1)] \x in {2, 3, 4, 5, 11, 12, 13, 14, 15, 17, 18, 19, 21, 23, 24, 25, 26, 27, 28, 29, 30}
\draw (v\y) to (v\x);
\draw (v5) to (v1);
\foreach \x in {6, 7, 8, 9, 10}
\draw (v5) to (v\x);
\draw (v6) to (v7);
\draw (v8) to (v9);
\draw (v12) to (v16);
\draw (v12) to (v20);
\draw (v19) to (v15);
\draw (v21) to (v15);
\draw (v22) to (v15);
\draw (v30) to (v24);
\end{tikzpicture}
\hskip0.5in
\begin{tikzpicture}[auto, scale=0.75]
\tikzstyle{vertex}=[draw, circle, inner sep=0.8mm]
\node (u1) at (0,0) [vertex,   label=left: $u_1$] {};
\node (w1) at (1.25,0) [vertex, label=right: $w_1$] {};
\node (u2) at (0,-0.5) [vertex,   label=left: $u_2$] {};
\node (w2) at (1.5,-0.5) [vertex, label=right:$w_2$] {};
\node (u3) at (0,-1) [vertex,   label=left: $u_3$] {};
\node (w3) at (1.25,-1) [vertex, label=right:$w_3$] {};
\node (u4) at (0,-1.5) [vertex,  label=left: $u_4$] {};
\node (w4) at (1.5,-1.5) [vertex, label=right:$w_4$] {};
\draw (u1) to (w1);
\draw (u2)  to (w1);
\draw(u3) to (w1);
\draw (u3) to (w3);
\draw (u4) to (w4);
\draw (u4) to (w3);
\draw (u4) to (w2);
\draw (w2) to (w3);
\draw (w1) to (w2);
\draw (w1) to (w3);
\draw (u2) to (w2);
\draw(w4) to (w3);
\draw(w4) to (w2);
\end{tikzpicture}
\end{center}
\caption{Two graphs that have $\lambda^-(G)=0$: a doubly extended bouquet dipole graph on the left and a two-set perfectly orderable graph on the right.}
\label{fig:skewzero}
\end{figure}
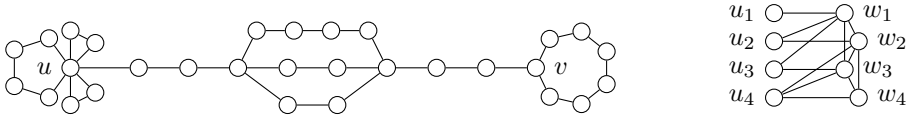

In addition, there exist graphs that have $\Zo(G)=0$ and therefore $\Fo(G)$ is undefined, specifically two-set perfectly orderable graphs.   
\begin{dfn}
We say that a graph $G$ is \emph{two-set perfectly orderable} if 
\begin{enumerate}
\item $V(G)$ can be partitioned into ordered sets, $U = \left\{ u_1, u_2, \ldots, u_m \right\}$ and $W = \left\{ w_1, w_2, \ldots, w_m \right\}$ such that $u_i v \in E(G)$ only if $v=w_j$ where $j \leq i$, and 
\item $u_i w_i \in E(G)$ for all $i \in \{1, 2, \ldots, m\}$.  
\end{enumerate}
\end{dfn}
This gives us the following characterization of graphs with $\lambda^-(G)=0$.  
\begin{thm}
$\lambda^-(G)=0$ if and only if  $G$ is one of the following graphs.
\begin{enumerate}
\item $K_n$ \label{item:kn}
\item $\overline{K_n}$ \label{item:isolated}
\item a doubly extended bouquet-dipole graph. \label{item:doubly} 
\item a collection of one or more odd cycles that intersect in exactly one vertex.  \label{item:odd}
\item a two-set perfectly orderable graph.   \label{item:perfectly}
\end{enumerate}
\end{thm}

\begin{proof}
In \cite{ansill2016failed} the graphs with $\Zo(G) < \Fo(G)$ were characterized and are precisely Graphs \ref{item:kn} through \ref{item:odd}.  Graphs with $\Zo(G)=0$ and $\Fo(G)$ undefined were characterized in the same paper as Graph \ref{item:perfectly}.  By Lemma \ref{lem:equivskew}, the results holds. \qed \end{proof}

Since we've seen that two-set perfectly orderable graphs not only have $\lambda^-(G)=0$, but they also are precisely the graphs with $\Zo(G)=0$, we have the following statement about their skew zero forcing polynomials.  

\begin{cor}
The only graphs with skew zero forcing polynomial $\mathcal{Z}^-(G;x)=\sum_{i=0}^n {n \choose i}x^i = (x+1)^n$ are two-set perfectly orderable graphs.  
\end{cor}

Finally, we characterize digraphs with $\lambda(D)=0$.  
\begin{thm}
A digraph $D$ has $\lambda(D)=0$ if and only if $D$ is one of the following.
\begin{enumerate}
\item a directed cycle. \label{characterizationdirectedcycle}
\item a regular tournament on 5 vertices. \label{tournament5}
\item a digraph obtained from $K_n$ by removing the arcs of  \label{characterizationremovecycles}
\begin{enumerate}
\item a collection of vertex-disjoint directed cycles each of length at least 3 that span $V$ ($n\geq 3$), \label{cyclesspan}
\item a collection of vertex-disjoint directed cycles each of length at least 3 that span $V \backslash \{v\}$ for some $v \in V$ ($n\geq 4$), or  \label{cyclesspanbutone}
\item $vu$ for some $u, v \in V$ and  a collection of vertex-disjoint directed cycles each of length at least 3 that span $V \backslash \{v\}$ ($n\geq 4$). \label{cyclesspanleaf}
\end{enumerate}
\item a digraph obtained from $K_{n-1}  \overrightarrow{\vee} \{v\}$ by removing the arcs of a collection of vertex-disjoint directed cycles  each of length at least 3 that span $K_{n-1}$ ($n \geq 4$).  \label{characterizationsinkcomplement}
\item  $K_j \overrightarrow{\vee} \overline{K_{\ell}}$ where $j \geq 2$ and $\ell \geq 0$. \label{characterizationcomplete}
\item $\overline{K_n}$. \label{characterizationisolated}
\end{enumerate}
\end{thm}

\begin{proof}
In \cite{adams2021failed}, the list of graphs in the statement of this theorem were shown to be the graphs that have $\Z(D) < \F(D)$.  We then apply Lemma \ref{lem:equivstandard}.\qed
\end{proof}


\subsubsection{Comments on graphs with standard zero forcing span of 1}
We provide here a list of graphs that have $\lambda(G)=1$.  We have neither a characterization, nor any reason to believe the list of graphs is complete.  First, we note an immediate but notable property of graphs with $\lambda(G)=1$, and introduce a few definitions.  
\begin{obs}
$\lambda(G)=1$ if and only if $\Z(G)=\F(G)$.  
\end{obs}
Recall that  a \emph{module} $S \subseteq V(G)$ is a set of vertices such that for any vertex $u \notin S$, either $uv \in E(G)$ for all $v \in S$, or $uv \notin E(G)$ for any $v \in S$. In \cite{fetcie2015failed} it was shown that $\F(G)=n-2$ if and only if $G$ has a module of order 2.
The \emph{union} of graphs $G$ and $H$, denoted $G \cup H$, is the graph with vertex set $V(G) \cup V(H)$ and edge set $E(G) \cup E(H)$.  The \emph{join} $G \vee H$ is simply $G \cup H$ with the addition of an edge between $v_G$ and $v_H$ for every $v_G \in V(G)$ and every $v_H \in V(H)$. 
\begin{prop}
The following graphs have $\lambda(G)=1$.
\begin{enumerate}
\item a complete multipartite graph $K_{n_1, n_2, \ldots n_k}$ where $n_1 \geq 2$. \label{item:multipart}
\item $K_m \cup K_n$ where $m, n \geq 2$. \label{item:union}
\item $K_n \backslash M$ where $n \geq 3$ and $M$ is a nonempty matching.  \label{item:matching}
\item $K_n \vee \overline{K_m}$ where $m \geq 2$.  \label{item:completesplit}
\item a path on $4$ vertices. \label{item:pathcycle}
\end{enumerate} \label{prop:span1}
\end{prop}

\begin{proof}
For the complete multipartite graph, Item \ref{item:multipart}, note that any pair of vertices in a single partite set forms a module of order $2$, so $F(G)=n-2$. For any $S\subseteq V(G)$ with $|S|\leq n-3$, $S$ is a failed zero forcing set since if two of the vertices in $V\backslash S$ are in a single partite set, then $S$ is a failed zero forcing set, and otherwise, each vertex of $V \backslash S$ is in a distinct partite set, and each vertex in $S$ is adjacent to at least two vertices in $V \backslash S$.  Take $u \in V(P_1)$ where $P_1$ is a partite set with at least two vertices, and $v \in P_2$ where $P_2$ is any other partite set.  Then $V \backslash \{u, v\}$ is a zero forcing set, giving us that $\F(G)=\Z(G)=n-2$, and $\lambda(G)=1$.  

For $K_m \cup K_n$, Item \ref{item:union}, note that any pair of vertices in the same component form a module of order $2$.  Thus, $F(G) = n+m-2$.  For $S$ to be a zero forcing set, it must contain $m-1$ vertices of $K_m$ and $n-1$ vertices of $K_n$, giving us that $\Z(G)=n+m-2$ as well.  

For Item \ref{item:matching}, $G=K_n \backslash M$ where $M$ is a matching, note that the endpoints of any edge in the matching form a module of order $2$, so $\F(G)=n-2$.  However, $\Z(G)=n-2$ as well, since taking $S= V\backslash \{u, v\}$ for any $u, v \in V$ where $uv$ is not an edge in the matching forms a zero forcing set, and any set of $S'$ with $|S'| \leq n-3$ is not a zero forcing set since for any $v \in S'$, $v$ is adjacent to at least two  vertices in $V \backslash S'$.

For Item \ref{item:completesplit}, note that any pair of vertices in either $K_n$ or $\overline{K_m}$ forms a module of order two, giving us $\F(G)=n-2$.  Since any set $S \subseteq V(  K_n \vee \overline{K_m})$ with $|S| \leq n-3$ has at least two vertices missing from $K_n$ or $\overline{K_m}$, we have that $\Z(G) \geq n-2$.  By picking a set $S'$ with $S'= V(  K_n \vee \overline{K_m}) \backslash \{u,v\}$ where $u \in K_n$ and $v \in \overline{K_m}$ we see $\Z(G)=n-2=\F(G)$.  

For Item \ref{item:pathcycle}, note that $\Z(P_4)=\F(P_4)=1$.  
\qed
\end{proof}

\subsubsection{Characterizations of graphs with high zero forcing spans}

We now characterize graphs and digraphs with high zero forcing spans.  Specifically, we characterize graphs that have $\lambda(G) \geq n-3$ and  digraphs that have $\lambda(D) \geq n-2$.  For skew zero forcing, we show that $\lambda^-(G) \neq n$ for any graph $G$ and characterize graphs that have $\lambda^-(G)=n-1$.  


\begin{thm} \label{prop:high}
Graphs with $\lambda(G)=n-1$ or $\lambda(G)=n-2$ can be characterized as follows.  
$$\lambda(G)=
\begin{cases}
n-1 \mbox{ if and only if } G=K_1 \\
n-2 \mbox{ if and only if } G = P_{n-1} \cup K_1 \mbox{ or }  G=P_3 
\end{cases}$$
\end{thm}

\begin{proof}
For $\lambda(G)=n-1$, we must have $\F(G)=n-1$ and $\Z(G)=1$.  From \cite{fetcie2015failed}, $\F(G)=n-1$ gives us that $G$ has an isolated vertex.  It is well known that $\Z(G)=1$ if and only if $G$ is a path.   Thus, $G=K_1$.  

If $\lambda(G)=n-2$, then either $\F(G)=n-1$ and $\Z(G)=2$, or $\F(G)=n-2$ and $\Z(G)=1$.  In the first case, $\F(G)=n-1$ implies that $G$ has an isolated vertex.  Since $G$ has an isolated vertex with $\Z(G)=2$, we must have that the other component of $G$ is a path, giving us a path and a single isolated vertex.  Note that we can construct a failed zero forcing set of $G$ with $n-1$ vertices by taking all vertices but the isolated vertex, and a zero forcing set with $2$ vertices by taking one end vertex of the path along with the isolated vertex.  

If $\F(G)=n-2$ and $\Z(G)=1$, we have that $G$ must be a path, but from \cite{fetcie2015failed} that there are two pendant vertices since $G$ is a tree with $\F(G)=n-2$, giving us that $G=P_3$. Note the middle vertex forms a failed zero forcing set of maximum order, and the end vertex a zero forcing set of minimum order.  \qed
\end{proof}

We pause here to recall definitions that are essential in some of the characterizations below.

\begin{dfn} 
We  say that $G$ is a graph of \emph{two parallel paths} if $G$ itself is not a path, and $V(G)$ can be partitioned into subsets $V_1$ and $V_2$ such that the subgraphs induced by $V_1$ and $V_2$ are paths, and $G$ can be drawn in the plane so that the paths induced by $V_1$ and $V_2$ are parallel line segments, and edges between $V_1$ and $V_2$ can be drawn as straight line segments  that do not cross.  Such a drawing is known as a \emph{standard drawing}.  
\end{dfn} 

\begin{lem}
A graph $G$ that consists of two parallel paths has a module of order 2 if and only if any standard drawing of $G$ consists of one of the following: \label{lem:twoparallelmodule}
\begin{enumerate}
	\item $P_1 \cup P_2$ or $P_1 \cup P_1$. \label{isolated}
	\item $P_1=\{x\}$ and $P_k$ where $k \geq 3$, and \label{parallelp1}
	\begin{enumerate}
		\item for some $u, v, w$ that form a subpath of $P_k$, $N(x)=\{u, w\}$ or $N(x)=\{u, v, w\}$,  or 
		\item for an end vertex $v$ of $P_k$ with neighbor $w$, $N(x)=\{w\}$ or $N(x)=\{v, w\}$.  
	\end{enumerate}
	\item $P_2=uv$ and $P_k$ with $N(u) \cap V(P_k) =N(v) \cap V(P_k)$.  \label{parallelp2}
	\item $P_3=uvw$ and $P_k$ where $N(u)=N(w)=\{v\}$ or $N(u)=N(w)=\{v, x\}$ for some $x$ on $P_k$.  \label{parallelp3}
	\item $P_k = v_1v_2v_3v_4 \cdots v_k$ and $P_j = w_1w_2w_3\cdots w_j$ where $k, j \geq 2$, and the edges between $P_k$ and $P_j$ are one of $\left\{ v_{k-1}w_1, v_k w_2 \right\}$, 
	$\left\{ v_{k-1}w_1, v_k w_2, v_{k-1}w_2\right\}$, or
	$\left\{ v_{k-1}w_1, v_k w_2,v_{k}w_1\right\}$. \label{parallelpk}  
\end{enumerate}
\end{lem}

\begin{proof}
Note that in Item \ref{isolated}, if $G=P_1 \cup P_2$, then $V(P_2)$ forms a module of order 2; if $G=P_1 \cup P_1$ then $V(G)$ itself is a module of order 2.     In Item \ref{parallelp1}, $\{v,x\}$ forms a module of order 2.  For Item \ref{parallelp2}, $\{u,v\}$ forms a module of order 2.  For Item \ref{parallelp3}, $\{u, w\}$ forms a module of order 2.  For Item \ref{parallelpk}, $\{v_k, w_1\}$ forms a module of order 2.  

Now assume that $G$ is two parallel paths and has a module of order 2.  We will show that $G$ is one of the graphs described in Items \ref{isolated} through \ref{parallelpk}. Consider a standard drawing of $G$.  Call the two paths $P_k$ and $P_j$. 

First, assume that $k,j \geq 4$, and that $\{u, v\}$ is a module of order 2.  Note that we cannot have that $u, v \in V(P_k)$ (without loss of generality) because then $u, v$ will have different neighbors along $P_k$.  Hence, we must have that $u \in V(P_k)$ and $v \in V(P_j)$.  If  $u$ is an interval vertex in $P_k$, then $v$ is adjacent to both vertices that are adjacent to $u$ along $P_k$, and there is also an edge between $u$ and the vertex (or vertices) adjacent to $v$; this edge will cross one of the edges from $v$ to the neighbors of $u$ which contradicts the definition of parallel paths.  Thus we must have that $u$ is an end vertex of $P_k$ and $v$ is an end vertex of $P_j$.

Recall  that we're considering a standard drawing, $P_k = v_1v_2v_3v_4 \cdots v_k$ and $P_j = w_1w_2w_3\cdots w_j$.  Note that we cannot have that $u= v_1$ and $v = w_1$, or $u=v_k$ and $v=w_j$, since then the edge from $u$ to the neighbor of $v$ will cross the edge from $v$ to the neighbor of $u$.  Thus, without loss of generality, we have $u=v_k$ and $v=w_1$.  To satisfy $\{ u, v\}$ being a module of order $2$, we then have that $v_{k-1} v \in E(G)$ and $w_2 u \in E(G)$.  That is, $v_{k-1}, v, w_2, u$ form a $C_4$.  Note then that we may have $uv \in E(G)$ or $v_{k-1} w_2 \in E(G)$, but not both, since the edges would cross, giving us Item \ref{parallelpk}.  

We now consider $j = 3$, so $P_j= uvw$.  Note that along $P_j$, $N(u)=N(w)=\{v\}$.  For $\{u, w\}$ to form a module of order 2, we must have that $u$ and $w$ have the same neighbors in $P_k$ as well.  Note that if $u$ and $w$ have more than one neighbor in $P_k$, then we will have crossed edges between $P_k$ and $P_j$.  Hence, $u$ and $w$ have at most one neighbor in $P_k$, giving us Item \ref{parallelp3}.  By the same arguments we made for the case when $j \geq 4$, the only other possibility for $j=3$ is if the graph satisfies Item \ref{parallelpk}.   

If $j=2$, let $P_2=uv$. For the case $k=1$, let $V(P_1)=\{x\}$.  Note either $ux, vx \in E(G)$ or $ux, vx \notin E(G)$ satisfying Item \ref{isolated} or \ref{parallelp2}.  
If $k \geq 2$, note that for $\{u, v\}$ to be a module of order $2$, then they must have the same neighborhood in $V(P_k)$, satisfying Item \ref{parallelp2}.  Note that by the definition of parallel paths, $|N(u)\cap V(P_k)| = |N(v) \cap V(P_k)| \in \{0, 1\}$, else edges from $u$ and $v$ to their neighbors in $P_k$ would cross.  If $u,v$ do not form a module, then without loss of generality $\{u, w\}$ form a module of order $2$ for some $w \in V(P_j)$. By the same arguments for the cases $k, j \geq 4$, $G$ must satisfy Item \ref{parallelpk}.  

Finally, suppose $j=1$.  If $k=1$, we have Item \ref{isolated}.  If $k=2$, we have the same situation just described for $j=2$ and $k=1$.  If $k =3$, we must have Item \ref{parallelp3} by the arguments for the case $j=3$.  If $k=4$, note that no two vertices of $P_k$ can form a module of order 2.  Thus we must have that $\{x, v\}$ form a module of order 2 where $\{x\}=V(P_j)$ and $v \in V(P_k)$.  Then $x$ must be adjacent to $N(v)$, and may be adjacent to $v$ as well, giving us Item  \ref{parallelp1}.  \qed
\end{proof}

If $S \subseteq V(G)$, then we denote the subgraph of $G$ induced by $S$ by $G[S]$.  We now characterize graphs with zero forcing span of $n-3$.

\begin{thm} 
$\lambda(G)=n-3$ if and only if $G$ is one of the following graphs.
\begin{enumerate}
\item two parallel paths with an additional $K_1$ component. \label{k1comp}
\item $P_4$ \label{p4} or $P_5$.
\item two parallel paths such that any standard drawing has one of the following forms: \label{parallel}
	\begin{enumerate}
	\item $P_1=\{x\}$ and $P_k$ where $k \geq 3$, and 
	\begin{enumerate}
		\item for some $u, v, w$ that form a subpath of $P_k$, $N(x)=\{u, w\}$ or $N(x)=\{u, v, w\}$, or 
		\item for an end vertex $v$ with neighbor $w$, $N(x)=\{w\}$ or $N(x)=\{v, w\}$.  
	\end{enumerate}
	\item $P_2=uv$ and $P_k$ with $N(u) \cap V(P_k) =N(v) \cap V(P_k)$, and if $k=1$, then $\{ux, vx\} \in E(G)$ where $x=V(P_1)$.   
	\item $P_3=uvw$ and $P_k$ where $N(u)=N(w)=\{v\}$ or $N(u)=N(w)=\{v, x\}$ for some $x$ on $P_k$.  
	\item $P_k = v_1v_2v_3v_4 \cdots v_k$ and $P_j = w_1w_2w_3\cdots w_j$ where $k, j \geq 2$,  and the edges between $P_k$ and $P_j$ are one of $\left\{ v_{k-1}w_1, v_k w_2 \right\}$, 
	$\left\{ v_{k-1}w_1, v_k w_2, v_{k-1}w_2\right\}$, or
	$\left\{ v_{k-1}w_1, v_k w_2,v_{k}w_1\right\}$.   
	\end{enumerate}
\end{enumerate}
\label{prop:nminus3}
\end{thm}

\begin{proof}
Using Observation \ref{obs:formula}, we know that $\lambda(G)=n-3$ implies that (I) $\F(G)=n-1$ and $\Z(G)=3$, (II) $\F(G)=n-2$ and $\Z(G)=2$, or (III) $\F(G)=n-3$ and $\Z(G)=1$.  

For (I) we know that $\F(G)=n-1$ if and only if $G$ contains an isolated vertex \cite{fetcie2015failed}, giving us that (I) is satisfied if and only if $G$ has an isolated vertex $v$ and  $\Z(G[V(G)\backslash\{v\}])=2$.  By  \cite{ferrero2019rigid}, $\Z(G[V(G)\backslash\{v\}])=2$ if and only if $G[V(G)\backslash\{v\}]$ is two parallel paths.  Hence, we have that $\F(G)=n-1$ and $\Z(G)=3$ if and only if Item \ref{k1comp} holds.

For (II), we know that $\F(G)=n-2$ if and only if $G$ has a module of order 2 \cite{fetcie2015failed} and no isolated vertices, and $\Z(G)=2$ if and only if $G$ is two parallel paths.  Hence, (II) holds if and only if $G$ is two parallel paths with a module of order 2 and no isolated vertices.  By Lemma \ref{lem:twoparallelmodule}, then (II) holds if and only if $G$ satisfies Item \ref{parallel}.

For (III),  $\Z(G)=1$ if and only if $G$ is a path.  From \cite{fetcie2015failed}, $\F(P_n)=n-3$ if and only if $n=4$ or $n=5$.   Hence, we have that $\F(G)=n-3$ and $\Z(G)=1$ if and only if Item \ref{p4}  holds.\qed
\end{proof}

We now characterize graphs of highest possible skew zero forcing span, and improve the trivial bound from Observation \ref{obs:skew} to a tight bound.  
\begin{thm}
$\lambda^-(G)=n-1$ if and only if $G$ consists of a (possibly empty) two-set perfectly orderable graph and an isolated vertex.  Also, $\lambda^-(G) \leq n-1$.  \end{thm}
\begin{proof}
By Observation \ref{obs:formula}, $\lambda^-(G)=n-1$ implies that $\Fo(G)=n-1$ and $\Zo(G)=1$, or $\Fo(G)=n-2$ and $\Zo(G)=0$.  If $\Zo(G)=0$, then $\Fo(G)$ is undefined, meaning that the only feasible case is that $\Fo(G)=n-1$ and $\Zo(G)=1$.  From \cite{ansill2016failed}, $\Fo(G)=n-1$ if and only if $G$ has an isolated vertex, $v$.  Note that one possibility is that $V(G)=\{v\}$.  Otherwise, since $\Zo(G)=1$, and any zero forcing set must contain $v$, we have that $\Zo(H) =0$ where $H$ is the subgraph induced by $V(G) \backslash \{v\}$.   Thus, $H$ must be a two-set perfectly orderable graph.  

We noted above that $\lambda^-(G) \leq n$.  If $\lambda^-(G)=n$, then $\Fo(G)=n-1$ and $\Zo(G)=0$. This gives us that $G$ is a two-set perfectly orderable graph with an isolated vertex, which is a contradiction since no two-set perfectly orderable graph has an isolated vertex.  Hence $\lambda^-(G) \leq n-1$.  \qed
\end{proof}

Finally, we turn to high zero forcing spans of digraphs.  First, we recall some definitions.  

A \emph{source} $v \in V(D)$ has $N^-(v)=\varnothing$. A  path $P= (v_1, v_2, \ldots, v_k)$ is \emph{Hessenberg} if $E(D)$ does not contain any arc of the form $(v_i, v_j)$ with $j > i+1$.  A simple digraph $D$ is a \emph{digraph of two parallel Hessenberg paths} if $D$ is not a Hessenberg path,  $V(D) = \left\{ u_1, \ldots u_r, v_1, \ldots, v_s\right\}$ where $(u_1, \ldots, u_r)$ and $(v_1, \ldots v_s)$ are nonempty Hessenberg paths, and there do not exist, $i, j, k, \ell$ with $i <j$ and $k<\ell$ such that $\left\{ (u_k, v_j), (v_i, k_{\ell})\right\} \subseteq E(D)$.  That is, there are no pairs of forward crossing arcs. 

\begin{thm}
$\lambda(D)=n-1$ if and only if $D$ is a Hessenberg path $(v_1, v_2, \ldots, v_n)$ such that $v_1$ is a source.\end{thm}

\begin{proof}
From \cite{hogben2010minimum}, $\Z(D) = 1$ if and only if $D$ is a Hessenberg path.   
From \cite{adams2021failed}, $\F(D)=n-1$ if and only if $D$ has a source.  Thus, $\lambda(D)=n-1$ if and only if $D$ is a Hessenberg path with a source.  The only vertex in a Hessenberg path $(v_1, \ldots, v_n)$ that can be a source is $v_1$, since $v_{i-1} \in N^-(v_i)$ for any $v_i$ with $i >1$.  Taking any Hessenberg path such that $v_1$ has no in-neighbors gives us a Hessenberg path with a source, $v_1$.   \qed
\end{proof}

\begin{figure}[h!]
\begin{center}
 \begin{tikzpicture}[auto, scale=0.8]
\tikzstyle{vertex}=[draw, circle, inner sep=0.8mm]
\tikzset{edge/.style = {->,> = Latex}}
\tikzset{curvededge/.style = {->,> = Latex, bend right}}
\node (u) at (0,0) [vertex, label= below:$v_1$] {};
\node (v) at (1,0) [vertex,  fill=blue, label= below: $v_2$] {};
\node (w) at  (2,0) [vertex,  fill=blue, label= below: $v_3$] {};
\node (x) at (3,0) [vertex,  fill=blue, label= below: $v_4$] {};
\node (y) at (4,0) [vertex, label= below: $v_5$] {};
\draw [edge] (u) to (v);
\draw [edge] (v) to (w);
\draw [edge] (w) to (x);
\draw [edge] (x) to (y);
\draw [curvededge] (x) to (v);
\draw [curvededge] (y) to (x);
\draw [curvededge] (x) to (u);
\end{tikzpicture}
\end{center}
\caption{A digraph with $\lambda(D)=n-2$. \label{fig:highdigraph}}
\end{figure}
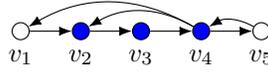

\begin{thm}
$\lambda(D)=n-2$ if and only if $D$ is one of
\begin{enumerate}
\item two parallel Hessenberg paths $(u_1, u_2, \ldots, u_k)$ and $(v_1, v_2, \ldots, v_j)$ such that $u_1$ or $v_1$ is a source, or \label{parallelHess}
\item a single Hessenberg path $(v_1, \ldots, v_n)$ with $N^-(v_1)=N^-(v_{i+1})=\{v_i\}$ for some $i$, $1<i<n$.  \label{twoneighborHess}
\end{enumerate}
\end{thm}

\begin{proof}
Note that from Observation \ref{obs:formula}, $\lambda(D)=n-2$ if and only if either $\Z(D)=1$ and $\F(D)=n-2$, or $\Z(D)=2$ and $\F(D)=n-1$.  From \cite{berliner2013minimum}, $\Z(D)=1$ if and only if $D$ is a Hessenberg path, and $\Z(D)=2$ if and only if $D$ consists of two parallel Hessenberg paths.  From \cite{adams2021failed}, $\F(D)=n-1$ if and only if $D$ has a source, and $\F(D)=n-2$ if and only if $D$ has no source and there exist $u, v \in V(D)$ with $N^-(u)\backslash\{v\} = N^-(v)\backslash\{u\}$.   

Noting that only the first vertex in a Hessenberg path can be a source, we have $\lambda(D)=n-2$ if and only if either Item \ref{parallelHess} holds, or $D$ is a Hessenberg path with no source and there exist $u, v \in V(D)$ with $N^-(u)\backslash\{v\} = N^-(v)\backslash\{u\}$.  By the definition of a Hessenberg path, $N^-(u)\backslash\{v\} = N^-(v)\backslash\{u\}$ is true if and only if  $u=v_1$ and $v= v_{i+1}$ where $1<i<n$, and $N^-(v_1) = N^-(v_{i+1}) = \{ v_i \}$, Item \ref{twoneighborHess}. \qed \end{proof}

The graph shown in Figure \ref{fig:highdigraph} is a Hessenberg path with $\lambda(D)=n-2$, since the blue vertices represent a failed zero forcing set $S$ with $|S|=n-2$, and $\{v_1\}$ forms a zero forcing set of order 1.  

\section{Zero forcing span characteristics for some graphs}

In this section, we look at the spans of trees, graphs with two or more components, and Cartesian products.  

\begin{prop}
For a tree $T$ on four or more vertices, $1 \leq \lambda(T) \leq n-3$ and these bounds are sharp.
\end{prop}
\begin{proof}
For any tree $T$, by Theorems \ref{prop:zero} and  \ref{prop:high}, $\lambda(T) \in \{0, n-2, n-1\}$ if and only if $T$ is $K_1, K_2$, or $P_3$.  Thus, if  $|V(T)| \geq 4$, we have $1 \leq \lambda(T) \leq n-3$.

From Proposition \ref{prop:span1}, for the star $T=K_{m,1}$ with $m \geq 2$, we see that $\lambda(K_{m,1}) = 1$, showing sharpness of the lower bound.  

From Theorem \ref{prop:nminus3}, if $T$ consists of a path on at least two vertices with two pendant vertices on one of its end vertices, we see that $\lambda(T)=n-3$.  Hence the bounds are sharp.  \qed
\end{proof}

For a disconnected graph, we display a formula for its span in terms of the zero forcing numbers, orders, and failed zero forcing numbers of its components.  
\begin{prop}
If $G$ is a graph with at least two components: $G_1, G_2, \ldots, G_k$, then 
$$\lambda(G)=|V(G)| +  \max_{1 \leq i \leq k} \left(  \F(G_i)  - |V(G_i)|\right) - \left(\sum_{i=1}^k \Z(G_i)\right) + 1$$
\end{prop}

\begin{proof}
From \cite{fetcie2015failed}, $\F(G) = |V(G)| +  \max_{1 \leq i \leq k} \left(  \F(G_i)  - |V(G_i)|\right)$.  For $S \subseteq V(G)$ to be a zero forcing set, $S \cap V(G_i)$ must be a zero forcing set for each $i$, $1 \leq i \leq k$, giving us that $\Z(G) =  \sum_{i=1}^k \Z(G_i)$.  Applying the formula from Observation \ref{obs:formula} completes the result. \qed
\end{proof}

We provide a bound on the zero forcing span of the Cartesian product of graphs $G$ and $H$, $G \square H$, in terms of the orders of $G$ and $H$ and their zero forcing and failed zero forcing numbers.  

\begin{prop}
Let $G \square H$ denote the Cartesian product of graphs $G$ and $H$.  Then
$$\lambda(G\square H)  \geq \max\{ \F(G) |V(H)|, F(H)|V(G)|\} - \min \{\Z(G) |V(H)|, \Z(H)|V(G)|\} + 1$$
\end{prop}

\begin{proof}
This follows from the bound on the failed zero forcing number of a Cartesian product \cite{fetcie2015failed} and the bound on the zero forcing number of a Cartesian product \cite{aim2008zero}, together with Observation \ref{obs:formula}. \qed
\end{proof}

\section{Conclusion}
In this paper, we introduced the idea of zero forcing span, and characterized graphs with high and low values in the context of standard zero forcing for both undirected and directed graphs, as well as in skew zero forcing.  Since the zero forcing span is intimately related to linear algebra and to zero forcing polynomials, further investigation of these relationships is a compelling direction, and could include further investigation of parameters studied here, or other variants including, for example, positive semidefinite zero forcing or power domination.  

\bibliographystyle{plain}

\end{document}